\documentclass[a4paper,11pt,twoside,reqno]{amsart}

\usepackage[utf8]{inputenc}
\usepackage[plainpages=false,pdfpagelabels=true]{hyperref}
\usepackage{amssymb,amsthm}
\usepackage[margin=1in]{geometry}
\usepackage{slashed}
\usepackage[dvipsnames]{xcolor}
\usepackage{comment}

\newtheorem{Satz}{Theorem}[section]
\newtheorem{Prop}[Satz]{Proposition}
\newtheorem{Lem}[Satz]{Lemma}

\newtheorem{Cor}[Satz]{Corollary}
\theoremstyle{definition}
\newtheorem{Dfn}[Satz]{Definition}
\newtheorem{Bem}[Satz]{Remark}
\newtheorem{Bsp}[Satz]{Example}
\newcommand{\tr}{\operatorname{Tr}}

\newcommand{\vol}{{\operatorname{vol}}}
\newcommand{\dv}{\text{ }dv}

\allowdisplaybreaks[1]

\renewcommand{\epsilon}{\varepsilon}

\newcommand{\R}{\ensuremath{\mathbb{R}}}
\newcommand{\N}{\ensuremath{\mathbb{N}}}

\newcommand{\s}{\ensuremath{\mathbb{S}}}
\numberwithin{equation}{section}

\usepackage{color}

\newtheorem{innercustomgeneric}{\customgenericname}
\providecommand{\customgenericname}{}
\newcommand{\newcustomtheorem}[2]{%
  \newenvironment{#1}[1]
  {%
   \renewcommand\customgenericname{#2}%
   \renewcommand\theinnercustomgeneric{##1}%
   \innercustomgeneric
  }
  {\endinnercustomgeneric}
}

\newcustomtheorem{customthm}{Theorem}

\title{A family of triharmonic maps to spheres in all dimensions greater than two}
\author{Volker Branding}
\date{\today}
\address{University of Vienna, Faculty of Mathematics\\
Oskar-Morgenstern-Platz 1, 1090 Vienna, Austria\\}
\email{volker.branding@univie.ac.at}

\author{Anna Siffert}
\address{Universität M\"unster, Mathematisches Institut\\
Einsteinstr. 62\\
48149 M\" unster\\
Germany}
\email{asiffert@uni-muenster.de}

\thanks{The first author gratefully acknowledges the support of the Austrian Science Fund (FWF) through the project "Geometric Analysis of Biwave Maps" (DOI: 10.55776/P34853)
}
\begin{document}

\begin{abstract}
We present a construction method for triharmonic maps to spheres.
In particular, we show that for any $m\in\mathbb{N}$ with $m\geq 3$ there exists a triharmonic map from $\mathbb{R}^m\setminus\{0\}$ into a round sphere. In addition, we provide a construction method for 
proper \(r\)-harmonic maps between spheres based on a suitable deformation of eigenmaps.
\end{abstract} 

\maketitle

\section{Introduction and Results}
Higher order variational problems are a research field of growing interest both in analysis and geometry. More often than not, it is a challenging task to develop a general existence theory
or to find classes of explicit solutions to geometric partial differential equations of arbitrary order.

\smallskip

The starting point of our study is comprised by the so-called \emph{harmonic maps}.
In order to define the latter we consider the energy for maps between two Riemannian manifolds $(M,g)$ and $(N,h)$ which is defined by
\begin{align}
\label{eq:energy}
E(\phi):=\frac{1}{2}\int_M|d\phi|^2\dv.
\end{align}
Here, we assume $M$ to be closed for simplicity of the presentation.
The critical points of \eqref{eq:energy} are characterized by the vanishing
of the \emph{tension field} which is given by
\begin{align}
\label{eq:tension}
\tau(\phi):=\tr\bar\nabla d\phi.
\end{align}
Here, \(\bar\nabla\) represents the connection on the pull-back bundle \(\phi^\ast TN\). Solutions of \(\tau(\phi)=0\) are precisely harmonic maps
and are characterized by a semilinear elliptic partial differential equation of second order. For an introduction to harmonic maps we refer to the book
\cite{MR2044031}.

\smallskip

A fourth order generalization of harmonic maps, which has received growing attention in recent years, is the concept of \emph{biharmonic maps}.
These are defined as critical points of the \emph{bienergy functional}, i.e.
\begin{align}
\label{eq:bienergy}
E_2(\phi):=\frac{1}{2}\int_M|\tau(\phi)|^2\dv.
\end{align}
More precisely, the critical points of \eqref{eq:bienergy} are those which satisfy
\begin{align}
\label{eq:bitension}
0=\tau_2(\phi):=\bar\Delta\tau(\phi)-\sum_{j=1}^mR^N(\tau(\phi),d\phi(e_j))d\phi(e_j),
\end{align}
where \(\{e_j\},j=1,\ldots,m\) represents a local orthonormal frame field tangent to \(M\) and \(\bar\Delta\) the connection Laplacian on the pull-back bundle \(\phi^\ast TN\). The quantity \(\tau_2(\phi)\) is called the
\emph{bitension field} of the map \(\phi\colon M\to N\).
It can be directly seen that every harmonic map gives a solution of the biharmonic map equation \eqref{eq:bitension}, hence one is interested in finding the non-harmonic solutions of \eqref{eq:bitension}
which are called \emph{proper biharmonic}.

The equation for biharmonic maps \eqref{eq:bitension} is a fourth order semilinear elliptic partial differential equation and due to the higher number of derivatives new technical difficulties in the mathematical analysis arise. For the current status of research on biharmonic maps we recommend the recent book \cite{MR4265170}.

\smallskip
A further variant of a higher order energy for maps between Riemannian manifolds is given by the trienergy functional which is defined as follows
\begin{align}
\label{eq:trienergy}
E_3(\phi)=\frac{1}{2}\int_M|\bar\nabla\tau(\phi)|^2\dv.
\end{align}

Critical points of \eqref{eq:trienergy} are called \emph{triharmonic maps} and are characterized by solutions of 
\begin{align}
\label{eq:triharmonic}\sum_{j=1}^m
\bar\Delta^2\tau(\phi)=\sum_{j=1}^m R^N(\bar\nabla_{e_j}\tau(\phi),\tau(\phi))d\phi(e_j)
-\sum_{j=1}^m R^N(\bar\Delta\tau(\phi),d\phi(e_j))d\phi(e_j).
\end{align}
Note that \eqref{eq:triharmonic} comprises a semilinear elliptic partial differential equation of order six. As in the case of biharmonic maps we are interested in finding those solutions of \eqref{eq:triharmonic} which are non-harmonic and refer to the latter as \emph{proper triharmonic}.

In order to develop the theory of triharmonic maps, it is first of all desirable to find explicit examples as a general existence theory seems to be out of reach.
So far, only few examples of triharmonic maps are known. In order to close this gap we present a construction method that 
allows us to manufacture countably many explicit triharmonic maps.

Below we focus on the case of spherical targets, for which equation \eqref{eq:triharmonic} acquires the form
\begin{align*}
\bar\Delta^2\tau(\phi)=&\sum_{j=1}^m\big(\langle\tau(\phi),d\phi(e_j)\rangle\bar\nabla_{e_j}\tau(\phi)-\langle\bar\nabla_{e_j}\tau(\phi),d\phi(e_j)\rangle\tau(\phi) \\
&-|d\phi|^2\bar\Delta\tau(\phi)+\langle d\phi(e_j),\bar\Delta\tau(\phi)\rangle d\phi(e_j)\big).
\end{align*}

For a classification of triharmonic maps from Euclidean space with finite energy we refer to \cite{MR4007262}.

\smallskip

Recently, the first author studied when maps of the form
\begin{align}
\label{maps}
    q=(\sin\alpha\cdot v,\cos\alpha),\qquad \alpha\in(0,\frac{\pi}{2}),
\end{align}
where $v:M\rightarrow\s^{n-1}\subset\mathbb{R}^n$
is a harmonic map, can be rendered biharmonic \cite{bra25}. 

In the major part of this manuscript we will consider a
\lq higher order analogue\rq\ of the radial projection map  which has been introduced in \cite{MR4371934,MR4593065} by Nakauchi.
In particular, in \cite{MR4593065}[Main Theorem, p.1] Nakauchi showed that for any $\ell,m\in\mathbb{N}$ with $\ell\leq m$ there exists a harmonic map
\begin{align*}
u^{(\ell)}\colon&\mathbb{R}^{m}\setminus\{0\}\rightarrow\s^{m^{\ell}-1},\\
\notag &x=(x_1,\dots,x_m)\mapsto u^{(\ell)}(x)=(u^{(\ell)}_{i_1\dots i_{\ell}}(x))_{1\leq i_1,\dots i_{\ell}\leq m}
\end{align*}
which we will define in detail in \eqref{nak-maps}.

Inspired by the approach of \cite{bra25} we will study maps of the form \eqref{maps} for both biharmonicity and triharmonicity.

Our first result in this regard is the following:
\begin{customthm}{\ref{thm:biharmonic-main}}
Let $m,\ell\in\mathbb{N}$.
The map
\(q:\mathbb{R}^{m}\setminus\{0\}\rightarrow \s^{m^{\ell}}\)
given by
\begin{align*}
q:=\big(\sin\alpha\cdot u^{(\ell)}_{i_1\dots i_{\ell}},\cos\alpha\big)
\end{align*}
is a proper biharmonic map if and only if
the following equation is satisfied 
\begin{align*}
\sin^2\alpha=\frac{\ell(\ell+m-2)+2m-8}{2\ell(\ell+m-2)}.
\end{align*}    
\end{customthm}

This allows us to obtain an existence result for biharmonic maps as follows: For each $m\in\mathbb{N}$ with $m\geq 3$ there exists a natural number $\ell\in\mathbb{N}$ with $\ell\leq m$ such that there is a proper biharmonic map $\mathbb{R}^{m}\setminus\{0\}\rightarrow \s^{m^{\ell}}$, see Corollary \ref{cor:biharmonic}.

\begin{Bem}
Theorem \ref{thm:biharmonic-main} was already known for \(\ell=1\),
see \cite{MR4076824}, and for \(\ell=2,3\), see \cite{MR4830603},
while Theorem \ref{thm:biharmonic-main} provides a generalization for all values of 
\(\ell\).
\end{Bem}

The second main contribution of this manuscript is a corresponding result for triharmonic maps.

\begin{customthm}{\ref{thm:triharmonic-main}}
Let $m,\ell\in\mathbb{N}$.
The map
\(q:\mathbb{R}^{m}\setminus\{0\}\rightarrow \s^{m^{\ell}}\)
given by
\begin{align*}
q:=\big(\sin\gamma\cdot u^{(\ell)}_{i_1\dots i_{\ell}},\cos\gamma\big)  
\end{align*}
is a triharmonic map if and only if the equation
\begin{align*}
&3\sin^4\gamma\,\ell^3(m+\ell-2)^3-2\sin^2\gamma\,\ell^2(m+\ell-2)^2[4+(m+\ell-6)(4+\ell)+\ell(m+\ell-2)] \\
\nonumber&+\ell(\ell+2)(\ell+4)(m+\ell-2)(m+\ell-4)(m+\ell-6)=0
\end{align*}
is satisfied.
\end{customthm}

As a consequence thereof we can deduce that
for each $m\in\mathbb{N}$ with $m\geq 3$ there exists a natural number $\ell\in\mathbb{N}$ with $\ell\leq m$ such that there is a triharmonic map $\mathbb{R}^{m}\setminus\{0\}\rightarrow \s^{m^{\ell}}$, see
Corollary \ref{cor:triharmonic} for the precise details.

\begin{Bem}
\begin{enumerate}
    \item Our result Theorem \ref{thm:triharmonic-main} requires that
the dimension of the domain is bigger than two. However, there are also existence results available for triharmonic maps from lower-dimensional domains.
For triharmonic curves on \(\s^n\), which correspond to triharmonic maps from a one-dimensional domain, we refer to \cite[Theorem 7]{MR4124860}, while further results on triharmonic surfaces can be found in \cite{MR4730425} and \cite{MR4462636}.
\item We could also consider the maps obtained in Theorem \ref{thm:triharmonic-main} as maps from the unit ball \(B^m\)
but in this case we would need to change to the \emph{weak} formulation
of the triharmonic map equation.
\end{enumerate}
\end{Bem}

Finally, we will study how the approach used in this manuscript
can be applied to a certain class of higher order energy functionals.

We define the \(r\)-energy \(E_r\) of a map between Riemannian manifolds as follows, where we need to distinguish between the even ($r=2s$) and the odd ($r=2s+1$) case:
\begin{align*}
E_{2s}(\phi)=&\int_M|\bar\Delta^{s-1}\tau(\phi)|^2\dv, \\
\nonumber E_{2s+1}(\phi)=&\int_M|\bar\nabla\bar\Delta^{s-1}\tau(\phi)|^2\dv.
\end{align*}
The Euler-Lagrange equations are given by a semilinear elliptic
partial differential equation of order \(2r\) and are presented in detail in Section 5. 
Again, it is obvious that a harmonic map represents a critical point of \(E_r\) and we thus say that a \(r\)-harmonic map is \emph{proper} if it is non-harmonic.
For the current status of research on higher order variational problems we refer to \cite{MR4106647}, a classification result
for polyharmonic maps from non-compact domains was achieved in 
\cite{MR4184658}.

Eventually we construct proper \(r\)-harmonic maps between spheres from eigenmaps.

\begin{customthm}{\ref{thm:r-harmonic-main}}
Let \(v\colon\s^m\to\s^{n-1}\) be an eigenmap, i.e. \(|\nabla v|^2=\lambda\) for some positive \(\lambda\in\R\).
Then, the map \(q:\s^m\to\s^n\) defined by
\begin{align*}
q:=(\sin\delta\cdot v,\cos\delta)
\end{align*}
is proper \(r\)-harmonic if and only if \(\sin\delta=\frac{1}{\sqrt{r}}\). These maps are unstable critical points of the \(r\)-energy \(E_r\).
\end{customthm}

\begin{Bem}
The proper \(r\)-harmonic maps given by the previous Theorem are 
very similar to the \(r\)-harmonic hyperspheres constructed in 
\cite[Theorem 1.1]{MR4106647}, \cite[Theorem 1.1]{MR3711937},
\cite[Example 6.2]{basi24} and also to the polyharmonic curve on the sphere provided by \cite[Theorem 1.5]{MR4542687}.
\end{Bem}

\medskip

\textbf{Sign conventions:}
Throughout this article we will employ the following sign conventions: 
For the Riemannian curvature tensor field we use 
$$
R(X,Y)Z=[\nabla_X,\nabla_Y]Z-\nabla_{[X,Y]}Z,
$$ 
where \(X,Y,Z\) are vector fields.

For the rough Laplacian on the pull-back bundle $\phi^{\ast} TN$ we employ the analysts sign convention, i.e.
$$
\bar\Delta = \tr(\bar\nabla\bar\nabla-\bar\nabla_\nabla).
$$
In particular, this implies that the Laplace operator has a negative spectrum.

\medskip

\textbf{Organization:}
In Section\,\ref{sec:prelim} we provide preliminaries.
We establish a construction method for proper biharmonic maps in Section\,\ref{sec:bih}. In Section\,\ref{sec:tri}, we generate triharmonic maps from biharmonic maps. Section\,\ref{sec:eigen}
contains a construction method for 
proper \(r\)-harmonic maps between spheres from eigenmaps.

\section{Preliminaries}
\label{sec:prelim}
The radial projection map is a harmonic map with a point singularity.
In \cite{MR4371934,MR4593065} \lq higher order analogues\rq\ of the radial projection map have been constructed. In this section we recall the definition of these maps. We further provide properties of these maps which will be used in Sections\,\ref{sec:bih} and \ref{sec:tri}.

\smallskip

In the case of maps to a spherical target, i.e. for maps $\phi:(M,g)\rightarrow (\s^n,g_R)$ with $g_R$ being the round metric, equation (\ref{eq:tension}) simplifies to
\begin{align}
\label{eq:harmonic-sphere}
\Delta \phi+|\nabla \phi|^2\phi=0,
\end{align}
where $\Delta$ denotes the Laplace-Beltrami operator.
Using equation \eqref{eq:harmonic-sphere} we will compute $\Delta^k \phi$, $k\in\mathbb{N}$, for specific maps $\phi$ which we introduce below.

\smallskip

We start by providing a preparatory result, which we will repeatedly make use of later and which can be deduced by a straightforward computation:

\begin{Lem}
\label{lem-1}
Let \(\Delta\) be the Laplacian on \(\R^m\) and \(r=\sqrt{x_1^2+\ldots x_m^2}\).
Then, the following identity holds
\begin{align}
\Delta\frac{1}{r^s}=\frac{s(s+2-m)}{r^{s+2}}    
\end{align}
for all \(s\in\N\).    
\end{Lem}

\medskip

We first focus on the case in which $\phi$ is the radial projection map $u$ on $\mathbb{R}^{m}\setminus\{0\}$,
i.e. we study the map
\(u\colon\mathbb{R}^{m}\setminus\{0\}\to\s^{m-1}\subset\R^{m}\)
given by
\begin{align}
\label{dfn:u}
u(x):=\frac{x}{r}, 
\end{align}
where $r:=r(x):=\lvert x\rvert$. In the next lemma we compute $\Delta^k u$.

\begin{Lem}
\label{lem:Delta-u}
Let \(u\colon\mathbb{R}^{m}\setminus\{0\}\to\s^{m-1}\subset\R^{m}\) be the map
defined in \eqref{dfn:u}. Then for each $k\in\mathbb{N}$ the following identity holds
\begin{align}
\label{eq:delta-k}
\Delta^k u&=\prod_{j=1}^k(2j-1-m)(2j-1)\frac{u}{r^{2k}}.
\end{align}
\end{Lem}
\begin{proof}
A direct calculation shows that $$\frac{\partial }{\partial x_j}\frac{x_i}{r}=\frac{\delta_{ij}}{r}-\frac{x_ix_j}{r^3}.$$ Hence Lemma\,\ref{lem-1} yields
\begin{align*}
\Delta u&=\frac{1-m}{r^2}u.
\end{align*}    
The claim then follows by induction, in which we make use of Lemma\,\ref{lem-1} and the identity
\begin{align}
\label{eq:nabla-nabla}
    \langle\nabla u,\nabla r^{-k}\rangle=0,
\end{align}
where $k\in\mathbb{N}$. 
\end{proof}

\medskip

We will now introduce the \lq higher order analogues\rq\ of the radial projection map  which has been provided in \cite{MR4371934,MR4593065}. 
In particular, Nakauchi \cite{MR4593065}[Main Theorem, p.1] showed that for any $\ell,m\in\mathbb{N}$ with $\ell\leq m$ there exists a harmonic map
\begin{align}
\label{nak-maps}
u^{(\ell)}\colon&\mathbb{R}^{m}\setminus\{0\}\rightarrow\s^{m^{\ell}-1},\\
\notag &x=(x_1,\dots,x_m)\mapsto u^{(\ell)}(x)=(u^{(\ell)}_{i_1\dots i_{\ell}}(x))_{1\leq i_1,\dots i_{\ell}\leq m}.
\end{align}
We set \(y_i=\frac{x_i}{r}\), then these maps can be defined recursively
\begin{align*}
u^{(1)}_{i_1}(x)=&y_{i_1},\\
u^{(\ell)}_{i_1\ldots i_\ell}(x)=&C_{\ell,m}\big(y_{i_1}u^{(\ell-1)}_{i_1\ldots i_{\ell-1}}(x)-\frac{1}{\ell+m-3}r\frac{\partial}{\partial x_{i_\ell}}u^{(\ell-1)}_{i_1\ldots i_{\ell-1}}(x)\big),
\end{align*}
where 
\begin{align*}
 C_{\ell,m}=\sqrt{\frac{\ell+m-3}{2\ell+m-4}}.   
\end{align*}
In particular, the maps $u^{(\ell)}_{i_1\dots i_{\ell}}$ have the following properties: 
\begin{enumerate}
\item $u^{(\ell)}_{i_1\dots i_{\ell}}$ satisfies   
\begin{align*}
    \Delta u^{(\ell)}_{i_1\dots i_{\ell}}(x)+\lvert\nabla u^{(\ell)}_{i_1\dots i_{\ell}}\rvert^2 u^{(\ell)}_{i_1\dots i_{\ell}}=0; 
\end{align*}
\item $u^{(\ell)}_{i_1\dots i_{\ell}}$ is a polynomial in $u_{i_1},\dots, u_{i_{\ell}}$ of degree $\ell$, where $u_{i_j}=\frac{x_{i_j}}{r}$;
\item $\lvert\nabla u^{(\ell)}_{i_1\dots i_{\ell}}\rvert^2=\frac{\ell(\ell+m-2)}{r^2}$.
\end{enumerate}
Furthermore, Nakauchi \cite{MR4593065}[Proposition 1, (1)] proved that these maps satisfy
\begin{align}
\label{ortho}
   \sum_{j=1}^m y_j\cdot\nabla_ju^{(\ell)}_{i_1\dots i_{\ell}}=0.
\end{align}

\begin{Lem}
\label{lem:delta-nak}
 For each $k\in\mathbb{N}$ the maps (\ref{nak-maps}) satisfy 
 \begin{align*}
          \Delta ^{k}u^{(\ell)}_{i_1\dots i_{\ell}}=\prod_{s=1}^k(2s+\ell-2)\prod_{j=1}^k(2j-\ell-m)\frac{u^{(\ell)}_{i_1\dots i_{\ell}}}{r^{2 k}}.
    \end{align*} 
\end{Lem}
\begin{proof}
We provide the claim by an induction argument.
For $k=1$ the claim is correct since $u^{(\ell)}_{i_1\dots i_{\ell}}$  satisfies properties $(1)$ and $(3)$ detailed above.

Now, assume that for a $k\in\mathbb{N}$ we have
    \begin{align*}
          \Delta ^{k}u^{(\ell)}_{i_1\dots i_{\ell}}=\prod_{s=1}^k(2s+\ell-2)\prod_{j=1}^k(2j-\ell-m)\frac{u^{(\ell)}_{i_1\dots i_{\ell}}}{r^{2k}}.
    \end{align*}
  With the help of \eqref{ortho} we can deduce
     \begin{align*}
          \Delta ^{k+1}u^{(\ell)}_{i_1\dots i_{\ell}}&=\prod_{s=1}^k(2s+\ell-2)\prod_{j=1}^k(2j-\ell-m)[\frac{1}{r^{2 k}}\Delta u^{(\ell)}_{i_1\dots i_{\ell}}+ \Delta(\frac{1}{r^{2 k}})u^{(\ell)}_{i_1\dots i_{\ell}}].
      \end{align*}    
       Using properties $(1)$ and $(3)$ of the maps (\ref{nak-maps}) as well as Lemma\,\ref{lem-1}, we hence obtain 
         \begin{align*}
         \Delta ^{k+1}u^{(\ell)}_{i_1\dots i_{\ell}} =&\prod_{s=1}^k(2s+\ell-2)\prod_{j=1}^k(2j-\ell-m)\\
         &\times
         [\frac{1}{r^{2k+2}}(2-m-\ell)\ell u^{(\ell)}_{i_1\dots i_{\ell}}+\frac{1}{r^{2k+2}}2k(2k+2-m)u^{(\ell)}_{i_1\dots i_{\ell}}]\\
          =&\prod_{s=1}^k(2s+\ell-2)\prod_{j=1}^k(2j-\ell-m)[2\ell-\ell^2-\ell m+4k^2+4k-2km]\frac{u^{(\ell)}_{i_1\dots i_{\ell}}}{r^{2k+2}}\\
          =&\prod_{s=1}^k(2s+\ell-2)\prod_{j=1}^k(2j-\ell-m)(2(k+1)+\ell-2)(2(k+1)-\ell-m)\frac{u^{(\ell)}_{i_1\dots i_{\ell}}}{r^{2 k+2}}\\
          =&\prod_{s=1}^{k+1}(2s+\ell-2)\prod_{j=1}^{k+1}(2j-\ell-m)\frac{u^{(\ell)}_{i_1\dots i_{\ell}}}{r^{2k+2}},
    \end{align*}
    which completes the induction argument.
\end{proof}

\begin{Bem}
Note that for $\ell=1$ Lemma\,\ref{lem:delta-nak} reduces 
to Lemma\,\ref{lem:Delta-u}.
\end{Bem}

Let us we briefly recall a number of facts on eigenmaps, for more details 
on this subject we refer to the books of Eells \(\&\) Ratto \cite[Part 3]{MR1242555} and Baird \(\&\) Wood \cite[p. 79]{MR2044031}.

\begin{Dfn}
A smooth map \(\phi\colon M\to\s^n\) is called eigenmap if the components
of \(u:=\iota\circ\phi\colon M\to\s^n\subset\R^{n+1}\) are all eigenfunctions of the Laplace-Beltrami operator on 
\(M\) with the same eigenvalue.
\end{Dfn}

A direct consequence of the previous definition is the next result:
\begin{Cor}
A smooth map \(\phi\colon M\to\s^n\) is an eigenmap if and only if it is harmonic
with constant energy density \(e(\phi):=\frac{1}{2}|d\phi|^2\).
\end{Cor}

\begin{Prop}
\label{prop:homo-pol}
Let \(F\colon\R^{m+1}\to\R^{n+1}\) be a harmonic map whose components are homogeneous
polynomials of the same degree \(k\in\{1,2,\ldots\}\).
Moreover, suppose that \(F\) restricts to a map \(\phi\colon\s^m\to\s^n\).
Then \(\phi\) is harmonic with constant norm of the gradient, i.e.
\begin{align}
|d\phi|^2=k(k+m-1).
\end{align}
\end{Prop}

\section{Biharmonic maps from harmonic maps}
\label{sec:bih}
In this section we present a method for constructing biharmonic maps from harmonic maps to spheres.

\smallskip

We consider the canonical embedding 
\(\iota\colon\s^n\hookrightarrow\R^{n+1}\) and will still write
\(u=(u_1,\ldots,u_{n+1})\) for \(\iota\circ u\).
In the following, we will employ the notation
\begin{align*}
\nabla u=(\nabla u_1,\ldots,\nabla u_{n+1}),
\qquad \Delta u=(\Delta u_1,\ldots,\Delta u_{n+1}).
\end{align*}
Here, \(\nabla\) refers to the gradient of \((M,g)\) and \(\Delta\)
is the Laplace-Beltrami operator. Note that each component of \(\nabla u\) is an \(m\)-dimensional vector, where \(m=\dim M\).

Recall that for a map \(u\colon M\to\s^{n}\) the bienergy functional
acquires the form
\begin{align*}
E_2(u)=\frac{1}{2}\int_M\big(|\Delta u|^2-|\nabla u|^4\big)\dv
\end{align*}
and a map is biharmonic if it satisfies
\begin{align}
\label{eq:biharmonic-sphere}
\Delta^2u+2\operatorname{div}\big(|\nabla u|^2\nabla u\big)
-\big(\langle\Delta^2u,u\rangle-2|\nabla u|^4\big)u=0.
\end{align}

In order to construct solutions of \eqref{eq:biharmonic-sphere} we now make the following ansatz: We consider the map \(q\colon M\to\s^{n}\subset\R^{n+1}\) 
given by
\begin{align}
\label{eq:definition-q-biharmonic}
q:=\big(\sin\alpha\cdot v,\cos\alpha\big),   
\end{align}
where \(v\colon M\to\s^{n-1}\subset\R^{n}\) is a given harmonic map, i.e. a solution of \eqref{eq:harmonic-sphere},
and \(\alpha\in (0,\frac{\pi}{2})\) is a real parameter.

The next proposition, which was already established in \cite[Proposition 3.1]{bra25}, gives a characterization when maps of the form \eqref{eq:definition-q-biharmonic} are actually biharmonic.

\begin{Prop}
The map \(q\colon M\to\s^n\subset\R^{n+1}\) defined in \eqref{eq:definition-q-biharmonic} is biharmonic if and only if the following equation holds
\begin{align}
\label{eq:biharmonic-q}
(\Delta|\nabla v|^2)v=-2\cos^2\alpha\,\nabla(|\nabla v|^2)\nabla v
+(1-2\sin^2\alpha)|\nabla v|^4v.
\end{align}
\end{Prop}

\begin{proof}
This follows from inserting the ansatz \eqref{eq:definition-q-biharmonic} into the equation
for biharmonic maps to spheres \eqref{eq:biharmonic-sphere} and using that \(v\colon M\to\s^{n-1}\) is a harmonic map. In particular, this implies that 
\begin{align*}
0=&\Delta^2v+\Delta|\nabla v|^2v+2\nabla|\nabla v|^2\nabla v
-|\nabla v|^4v,    
\end{align*}
see \cite[Lemma 2.1]{bra25}, which already completes the proof.
\end{proof}

\begin{Satz}
\label{thm:biharmonic-main}
Let $m,\ell\in\mathbb{N}$.
The map
\(q:\mathbb{R}^{m}\setminus\{0\}\rightarrow \s^{m^{\ell}}\)
given by
\begin{align*}
q:=\big(\sin\alpha\cdot u^{(\ell)}_{i_1\dots i_{\ell}},\cos\alpha\big),
\end{align*}
where \(u^{(\ell)}_{i_1\dots i_{\ell}}\) is defined in \eqref{nak-maps},
is a proper biharmonic map if and only if
the following equation is satisfied 
\begin{align}
\label{eq:constraint-biharmonic}
\sin^2\alpha=\frac{\ell(\ell+m-2)+2m-8}{2\ell(\ell+m-2)}.
\end{align}    
\end{Satz}

\begin{proof}
We have 
\begin{align*}
|\nabla u^{(\ell)}|^2=\frac{\ell(\ell+m-2)}{r^2}
\end{align*}
and the result follows from a direct calculation using \eqref{eq:biharmonic-q}.
\end{proof}

\begin{Bem}
The proper biharmonic maps provided by Theorem \ref{thm:biharmonic-main} were already known for \(\ell=1\) \cite[Theorem 1.1]{MR4076824} and \(\ell=2,3\) \cite[Theorems 1.2, 1.4]{MR4830603}.
\end{Bem}

In the following lemma we discuss under which conditions on $m$ and $\ell$ identity \eqref{eq:constraint-biharmonic} admits solutions.
Note however that not every solution of this identity gives rise to a proper biharmonic map: In the construction of Nakauchi, which we recalled in Section\,\ref{sec:prelim}, we have the condition $\ell\leq m$. In other words, solutions of \eqref{eq:constraint-biharmonic} with $\ell>m$ do not give rise to proper biharmonic maps.

\begin{Lem}
\label{lem:bih}
For each $(\ell,m)\in\mathbb{N}\times\mathbb{N}$ satisfying one of the following conditions
\begin{enumerate}
\item $\ell=1$ and $m\in\{4,5,6\}$;
\item $\ell=2$ and $m\geq 3$;
    \item $\ell\geq 3$ and $m\geq 2$;
\item $\ell\geq 4$ and $m=1$ 
\end{enumerate}
identity \eqref{eq:constraint-biharmonic} admits a solution.  For all other combinations of $(\ell,m)\in\mathbb{N}\times\mathbb{N}$ identity \eqref{eq:constraint-biharmonic} does not admit a solution.  
\end{Lem}

\begin{proof}
We will conduct a case distinction as follows:
\begin{enumerate}
    \item Let $\ell=1$, then
\begin{align*}
 \frac{\ell(\ell+m-2)+2m-8}{2\ell(\ell+m-2)}=\frac{(m-1)+2m-8}{2(m-1)}=\frac{1}{2}+\frac{m-4}{m-1}.  
\end{align*}
Hence for $\ell=1$, \eqref{eq:constraint-biharmonic} admits solutions if and only if $m\in\{4,5,6\}$.

\item Let $\ell=2$, then
\begin{align*}
 \frac{\ell(\ell+m-2)+2m-8}{2\ell(\ell+m-2)}=\frac{m-2}{m}=1-\frac{2}{m}. 
\end{align*}
Hence for $\ell=2$, \eqref{eq:constraint-biharmonic} admits solutions if and only if $m\geq 3$.

\item For each $(\ell,m)$ with $\ell\geq 3$ and $m\geq 2$ we have 
\begin{align*}
 \frac{\ell(\ell+m-2)+2m-8}{2\ell(\ell+m-2)}=\frac{1}{2}+ \frac{m-4}{\ell(\ell+m-2)}\leq \frac{1}{2}+\frac{m-4}{3(m+1)}<1  
\end{align*}
and
\begin{align*}
 \frac{\ell(\ell+m-2)+2m-8}{2\ell(\ell+m-2)}=\frac{1}{2}+ \frac{m-4}{\ell(\ell+m-2)}\geq \frac{1}{2}-\frac{2}{\ell^2}>0.  
\end{align*}
Hence, assuming the constraints $\ell\geq 3$ and $m\geq 2$ equation \eqref{eq:constraint-biharmonic} admits a solution.

For $(\ell,m)=(3,1)$ we get
\begin{align*}
 \frac{\ell(\ell+m-2)+2m-8}{2\ell(\ell+m-2)}=0,  
\end{align*}
therefore in this case \eqref{eq:constraint-biharmonic} does not admit any solutions.

\item For each $(\ell,m)$ with $\ell\geq 4$ and $m=1$ we have 
\begin{align*}
 \frac{\ell(\ell+m-2)+2m-8}{2\ell(\ell+m-2)}=\frac{1}{2}-\frac{3}{\ell(\ell-1)}> 0 
\end{align*}
and
\begin{align*}
 \frac{\ell(\ell+m-2)+2m-8}{2\ell(\ell+m-2)}&=\frac{1}{2}-\frac{3}{\ell(\ell-1)}<1. 
\end{align*}
\end{enumerate}
\end{proof}

\begin{Bsp}
For $m=\ell=4$, \eqref{eq:constraint-biharmonic}
admits the solution $\alpha=\frac{\pi}{4}$.
Hence the map 
\(q:\mathbb{R}^{4}\setminus\{0\}\rightarrow \s^{256}\)
given by
\begin{align*}
q:=\big(\frac{1}{\sqrt{2}}\cdot u^{(4)},\frac{1}{\sqrt{2}}\big),   
\end{align*}
is proper biharmonic.
\end{Bsp}

From Theorem\,\ref{thm:biharmonic-main} and Lemma\,\ref{lem:bih} we obtain the following corollary.

\begin{Cor}
\label{cor:biharmonic}
For each $m\in\mathbb{N}$ with $m\geq 3$ there exists a natural number $\ell\in\mathbb{N}$ with $\ell\leq m$ such that there is a proper biharmonic map $\mathbb{R}^{m}\setminus\{0\}\rightarrow \s^{m^{\ell}}$.
\end{Cor}

\section{Triharmonic maps from harmonic maps}
\label{sec:tri}
As a first step towards deriving an explicit solution of the equation for triharmonic maps
we need to rewrite the trienergy \eqref{eq:trienergy} in the case that \(N=\s^n\subset\R^{n+1}\).
To this end consider \(\phi\colon M\to\s^n\), the canonical inclusion 
\(\iota\colon\s^n\hookrightarrow\R^{n+1}\) and the composite map \(u:=\iota\circ\phi\colon M\to\R^{n+1}\).
Then, we have the identities
\begin{align}
\label{eq:trienergy-embed}
d\iota(\tau(\phi))=&\Delta u +|\nabla u|^2u,\\
\nonumber d\iota(\bar\nabla\tau(\phi))=&\nabla\Delta u +\langle\Delta u,\nabla u\rangle u
+\nabla(|\nabla u|^2u),
\end{align}
which allow us to establish the following result:

\begin{Lem}
The trienergy of a smooth map \(u\colon M\to\s^n\subset\R^{n+1}\) is given by \begin{align}
 \label{eq:trienergy-extrinsic}
 E_3(u)=\frac{1}{2}\int_M\big(&|\nabla\Delta u|^2+|\langle\Delta u,\nabla u\rangle|^2
+\big|\nabla|\nabla u|^2\big|^2+|\nabla u|^6 \\
\nonumber&-2|\langle\nabla\Delta u,u\rangle|^2
-2|\Delta u|^2|\nabla u|^2
 \big)\dv.
 \end{align}
\end{Lem}
\begin{proof}
Using \eqref{eq:trienergy-embed} a direct calculation shows that
\begin{align}
\label{eq:tri}
|d\iota(\bar\nabla\tau(\phi))|^2=&|\nabla\Delta u|^2+|\langle\Delta u,\nabla u\rangle|^2 +|\nabla(|\nabla u|^2u)|^2 \\
\notag&+2\langle\nabla\Delta u,u\rangle\langle\Delta u,\nabla u\rangle
+2\langle\nabla\Delta u,\nabla(|\nabla u|^2u)\rangle
+2\langle\Delta u,\nabla u\rangle\langle u,\nabla(|\nabla u|^2u)\rangle.
\end{align}
The first two summands on the right hand side do not require any further manipulations in order to deduce \eqref{eq:trienergy-extrinsic}.
Due to the assumption of a spherical target, we directly obtain that the third summand on the right hand side of \eqref{eq:tri} is given by
\begin{align*}
|\nabla(|\nabla u|^2u)|^2=&\big|\nabla|\nabla u|^2\big|^2+|\nabla u|^6.
\end{align*}
Using integration by parts, the integral over $M$ of the fourth summand on the right hand side of \eqref{eq:tri} reads
\begin{align*}
\int_M \langle\nabla\Delta u,u\rangle\langle\Delta u,\nabla u\rangle\dv
=&-\int_M\big(|\langle\nabla\Delta u,u\rangle|^2
+\langle\Delta^2u,u\rangle\langle\Delta u,u\rangle
+\langle\nabla\Delta u,\nabla u\rangle\langle\Delta u,u\rangle\big)\dv \\
=&-\int_M|\langle\nabla\Delta u,u\rangle|^2\dv
+\int_M\langle\Delta^2u,u\rangle|\nabla u|^2\dv \\
&+\int_M\langle\nabla\Delta u,\nabla u\rangle|\nabla u|^2\dv.
\end{align*}
In addition, using integration by parts we find that the integral over $M$ of the fifth summand on the right hand side of \eqref{eq:tri} is given by
\begin{align*}
\int_M\langle\nabla\Delta u,\nabla(|\nabla u|^2u)\rangle \dv=&        
-\int_M|\nabla u|^2\langle\Delta^2u,u\rangle\dv.
\end{align*}
Using $\langle u,\nabla u\rangle=0$ and integration by parts once more, we find that the integral over $M$ of the sixth summand on the right hand side of \eqref{eq:tri} acquires the form
\begin{align*}
\int_M\langle\Delta u,\nabla u\rangle\langle u,\nabla(|\nabla u|^2u)\rangle \dv
=&\int_M\langle\Delta u,\nabla u\rangle \nabla|\nabla u|^2 \dv \\
=&-\int_M\langle\nabla\Delta u,\nabla u\rangle|\nabla u|^2\dv 
-\int_M|\Delta u|^2|\nabla u|^2\dv.
\end{align*}
The claim now follows by taking the integral over $M$ of the right hand side of \eqref{eq:tri} and plugging in the identities for the last four summands which we just derived.
\end{proof}

In the next proposition we characterize critical points of the trienergy
\eqref{eq:trienergy-extrinsic}.

\begin{Prop}
Let \(u\colon M\to\s^n\subset\R^{n+1}\) be a smooth map.
The critical points of the trienergy \eqref{eq:trienergy-extrinsic}    
are characterized by solutions of the identity
\begin{align}
\label{eq:triharmonic-u}
\Delta^3u=&\Delta\big(\langle\Delta u,\nabla u\rangle\nabla u\big)
-\nabla\big(\langle\Delta u,\nabla u\rangle\Delta u\big)
+2\nabla\big((\Delta|\nabla u|^2)\nabla u\big)-3\nabla\big(|\nabla u|^4\nabla u\big) \\
\nonumber&+2\Delta\nabla\big(\langle\nabla\Delta u,u\rangle u\big)
-2\langle\nabla\Delta u,u\rangle\nabla\Delta u
-2\Delta\big(\Delta u|\nabla u|^2\big)+2\nabla\big(|\Delta u|^2\nabla u\big) \\
\nonumber&+\big[
\langle\Delta^3u,u\rangle
-\langle u,\Delta\big(\langle\Delta u,\nabla u\rangle\nabla u\big)\rangle
+\langle u,\nabla\big(\langle\Delta u,\nabla u\rangle\Delta u\big)\rangle
+2|\nabla u|^2\Delta|\nabla u|^2 \\
\nonumber&-3|\nabla u|^6 
-2\langle u,\Delta\nabla\big(\langle\nabla\Delta u,u\rangle u\big)\rangle
+2|\langle u,\nabla\Delta u\rangle|^2 \\
\nonumber&+2\langle u,\Delta\big(\Delta u|\nabla u|^2\big)\rangle 
+2|\Delta u|^2|\nabla u|^2
\big]u.
\end{align}
\end{Prop}
\begin{proof}
We consider a smooth one-parameter variation of the map \(u\colon M\to\s^n\subset\R^{n+1}\), that is a smooth map \(u_t\colon (-\epsilon,\epsilon)\times M\to\s^n\subset\R^{n+1}\), where \(\epsilon>0\) and $u_0=u$.
The variation vector field of this variation is denoted by $\eta$, i.e.
\begin{align*}
\frac{du_t}{dt}\big|_{t=0}=\eta.    
\end{align*}
In order to take care of the fact that we are considering a spherical target
we introduce a Lagrange multiplyer \(\lambda\in\R\) into the trienergy \eqref{eq:trienergy-extrinsic}, i.e. we deal with the energy functional 
\begin{align}
\label{eq:lem-var-a}
E_3(u)=\frac{1}{2}\int_M\big(&|\nabla\Delta u|^2+|\langle\Delta u,\nabla u\rangle|^2
+\big|\nabla|\nabla u|^2\big|^2+|\nabla u|^6 \\
\nonumber&-2|\langle\nabla\Delta u,u\rangle|^2
-2|\Delta u|^2|\nabla u|^2 +\lambda(|u|^2-1)
 \big)\dv,    
\end{align}
which we also denote by $E_3$ for shortness of notation. We consider the variation of this functional and
calculate the variation of each of the summands in 
\eqref{eq:lem-var-a} separately.
First of all, it is straightforward to see that the variation of the first summand is given by
\begin{align*}
\frac{d}{dt}\big|_{t=0}\frac{1}{2}\int_M|\nabla\Delta u_t|^2\dv
=-\int_M\langle\eta,\Delta^3u\rangle\dv.    
\end{align*}
Furthermore, we calculate the variation of the second summand of \eqref{eq:trienergy-extrinsic} as follows
\begin{align*}
\frac{d}{dt}\big|_{t=0}\frac{1}{2}\int_M|\langle\Delta u_t,\nabla u_t\rangle|^2\dv 
&=\int_M\langle\Delta u,\nabla u\rangle
\big(\langle\Delta\eta,\nabla u\rangle+\langle\Delta u,\nabla\eta\rangle\big)\dv \\
&=\int_M\langle\Delta\big(\langle\Delta u,\nabla u\rangle\nabla u\big),\eta\rangle\dv 
-\int_M\langle\nabla\big(\langle\Delta u,\nabla u\rangle\Delta u\big),\eta\rangle\dv.
\end{align*}
The variation of the third summand of \eqref{eq:trienergy-extrinsic} is given by
\begin{align*}
\frac{d}{dt}\big|_{t=0}\frac{1}{2}\int_M\big|\nabla|\nabla u_t|^2\big|^2\dv
=&-2\int_M\Delta|\nabla u|^2\langle\nabla\eta,\nabla u\rangle\dv \\
=&2\int_M\langle\nabla(\Delta|\nabla u|^2\nabla u),\eta\rangle\dv.
\end{align*}
Moreover, we find that the variation of the fourth summand reads
\begin{align*}
\frac{d}{dt}\big|_{t=0}\frac{1}{2}\int_M|\nabla u_t|^6\dv=
-3\int_M\langle\nabla\big(|\nabla u|^4\nabla u\big),\eta\rangle\dv,    
\end{align*}
while the variation of the fifth summand yields
\begin{align*}
\frac{d}{dt}\big|_{t=0}\frac{1}{2}\int_M|\langle\nabla\Delta u_t,u_t\rangle|^2\dv
=&\int_M\langle\nabla\Delta u,u\rangle
\big(\langle\nabla\Delta\eta,u\rangle+\langle\nabla\Delta u,\eta\rangle\big)\dv \\
=&\int_M\big\langle-\Delta\nabla\big(\langle\nabla\Delta u,u\rangle u\big)
+\langle\nabla\Delta u,u\rangle\nabla\Delta u,\eta\big\rangle\dv.
\end{align*}
Finally, we obtain
\begin{align*}
\frac{d}{dt}\big|_{t=0}\frac{1}{2}\int_M|\Delta u_t|^2|\nabla u_t|^2\dv
&=\int_M\Delta u\Delta\eta |\nabla u|^2\dv+\int_M|\Delta u|^2\langle\nabla u,\nabla\eta\rangle \dv\\
&=\int_M\big\langle\Delta\big(\Delta u|\nabla u|^2\big)-\nabla\big(|\Delta u|^2\nabla u\big),\eta\big\rangle\dv.
\end{align*}
Collecting all the above terms the first variation of \eqref{eq:lem-var-a}
leads to the system
\begin{align*}
\Delta^3u=&\Delta\big(\langle\Delta u,\nabla u\rangle\nabla u\big)
-\nabla\big(\langle\Delta u,\nabla u\rangle\Delta u\big)
+2\nabla\big((\Delta|\nabla u|^2)\nabla u\big)-3\nabla\big(|\nabla u|^4\nabla u\big) \\
&+2\Delta\nabla\big(\langle\nabla\Delta u,u\rangle u\big)
-2\langle\nabla\Delta u,u\rangle\nabla\Delta u
-2\Delta\big(\Delta u|\nabla u|^2\big)+2\nabla\big(|\Delta u|^2\nabla u\big)+\lambda u,\\
|u|^2=&1.
\end{align*}
In order to eliminate \(\lambda\) from the first equation we take the scalar product with \(u\) and find
\begin{align*}
\lambda=&\langle\Delta^3u,u\rangle
-\langle u,\Delta\big(\langle\Delta u,\nabla u\rangle\nabla u\big)\rangle
+\langle u,\nabla\big(\langle\Delta u,\nabla u\rangle\Delta u\big)\rangle
+2|\nabla u|^2\Delta|\nabla u|^2 \\
&-3|\nabla u|^6 
-2\langle u,\Delta\nabla\big(\langle\nabla\Delta u,u\rangle u\big)\rangle
+2|\langle u,\nabla\Delta u\rangle|^2 \\
&+2\langle u,\Delta\big(\Delta u|\nabla u|^2\big)\rangle 
+2|\Delta u|^2|\nabla u|^2.
\end{align*}
Reinserting this choice of \(\lambda\) into the above system completes the proof.
\end{proof}

There are of course many equivalent ways in which the equation for triharmonic maps can be expressed. In the following lemma we slightly rewrite \eqref{eq:triharmonic-u} in a way best suited for our purposes.

\begin{Bem}
The version of the triharmonic map equation presented in 
\eqref{eq:triharmonic-u} may not be optimal for analytic studies of the subject
as it contains derivatives of order six on the right hand side.
However, using the geometry of the sphere (i.e. by differentiating \(|u|^2=1\)),
we can derive the following two equations
\begin{align*}
\langle\nabla\Delta u,u\rangle &
=-\nabla|\nabla u|^2-\langle\nabla u,\Delta u\rangle   
,\\
\langle\Delta^3u,u\rangle&=-\Delta^2|\nabla u|^2-\Delta|\Delta u|^2
-2\Delta\langle\nabla u,\nabla\Delta u\rangle
-\langle\Delta u,\Delta^2 u\rangle-2\langle\nabla u,\nabla\Delta^2u\rangle.
\end{align*}    
Applying these identities it is then possible to turn \eqref{eq:triharmonic-u}
into a semilinear equation.
\end{Bem}

\begin{Lem}
Let \(\tilde u\colon M\to\s^{n}\subset\R^{n+1}\) be a map of the form
\begin{align}
\label{dfn:tilde-u}
\tilde u:=(\sin\alpha\cdot u,\cos\alpha),\qquad \alpha\in(0,\frac{\pi}{2}),    
\end{align}
where \(u\colon M\to\s^{n-1}\subset\R^{n}\) is a given smooth map.
Then \(\tilde u\) is triharmonic if and only if \(\tilde u\) is a solution of
\begin{align}
\label{eq:triharmonic-tilde-u}
\Delta^3\tilde u=&\Delta\big(\langle\Delta\tilde u,\nabla\tilde u\rangle\nabla\tilde u\big)
-\nabla\big(\langle\Delta\tilde u,\nabla\tilde u\rangle\Delta\tilde  u\big) 
+2\nabla\big((\Delta|\nabla\tilde u|^2)\nabla\tilde u\big)
-3\nabla\big(|\nabla\tilde u|^4\nabla\tilde u\big) \\
\nonumber&+4(\nabla\langle\Delta^2\tilde u,\tilde u\rangle)\nabla\tilde u
+4\nabla(\langle\nabla\Delta \tilde u,\nabla \tilde u\rangle)\nabla \tilde u 
+2(\nabla\langle\nabla\Delta\tilde u,\tilde u\rangle)\Delta\tilde u
+2\Delta(\langle\nabla\Delta \tilde u,\tilde u\rangle\nabla\tilde u)
\\
&
\nonumber-2\langle\nabla\Delta\tilde u,\tilde u\rangle\nabla\Delta\tilde  u
-2\Delta\big(\Delta\tilde u|\nabla\tilde u|^2\big)+2\nabla\big(|\Delta \tilde u|^2\nabla\tilde u\big).    
\end{align}
\end{Lem}

\begin{proof}
We employ the equation for triharmonic maps to spheres \eqref{eq:triharmonic-u}
derived in the previous proposition.
For this purpose we rewrite the fifth summand of the right hand side of \eqref{eq:triharmonic-u} as follows
\begin{align*}
\Delta\nabla\big(\langle\nabla\Delta u,u\rangle u\big)
=&(\Delta\nabla\langle\nabla\Delta u,u\rangle)u
+2(\nabla\langle\Delta^2u,u\rangle)\nabla u
+2\nabla(\langle\nabla\Delta u,\nabla u\rangle)\nabla u \\
&+(\nabla\langle\nabla\Delta u,u\rangle)\Delta u
+\Delta(\langle\nabla\Delta u,u\rangle\nabla u).
\end{align*}
Inserting this expression into \eqref{eq:triharmonic-u} yields
\begin{align*}
\Delta^3u=&\Delta\big(\langle\Delta u,\nabla u\rangle\nabla u\big)
-\nabla\big(\langle\Delta u,\nabla u\rangle\Delta u\big)
+2\nabla\big((\Delta|\nabla u|^2)\nabla u\big)
-3\nabla\big(|\nabla u|^4\nabla u\big) \\
\nonumber&
+4(\nabla\langle\Delta^2u,u\rangle)\nabla u
+4\nabla(\langle\nabla\Delta u,\nabla u\rangle)\nabla u \\
&+2(\nabla\langle\nabla\Delta u,u\rangle)\Delta u
+2\Delta(\langle\nabla\Delta u,u\rangle\nabla u)
\\
&
-2\langle\nabla\Delta u,u\rangle\nabla\Delta u
-2\Delta\big(\Delta u|\nabla u|^2\big)+2\nabla\big(|\Delta u|^2\nabla u\big) \\
\nonumber&+\big[
\langle\Delta^3u,u\rangle
-\langle u,\Delta\big(\langle\Delta u,\nabla u\rangle\nabla u\big)\rangle
+\langle u,\nabla\big(\langle\Delta u,\nabla u\rangle\Delta u\big)\rangle
+2|\nabla u|^2\Delta|\nabla u|^2 \\
\nonumber&-3|\nabla u|^6 
-2\langle u,\Delta\nabla\big(\langle\nabla\Delta u,u\rangle u\big)\rangle
+2|\langle u,\nabla\Delta u\rangle|^2 \\
\nonumber&+2\langle u,\Delta\big(\Delta u|\nabla u|^2\big)\rangle 
+2|\Delta u|^2|\nabla u|^2
+2\Delta\nabla\langle\nabla\Delta u,u\rangle
\big]u.    
\end{align*}
Hence, we can rewrite the equation for triharmonic maps to spheres as 
\begin{align*}
P_3(u)-\langle P_3(u),u\rangle u=0,    
\end{align*}
where 
\begin{align*}
P_3(u)=&  
-\Delta^3u+\Delta\big(\langle\Delta u,\nabla u\rangle\nabla u\big)
-\nabla\big(\langle\Delta u,\nabla u\rangle\Delta u\big) 
+2\nabla\big((\Delta|\nabla u|^2)\nabla u\big)
-3\nabla\big(|\nabla u|^4\nabla u\big) \\
&+4(\nabla\langle\Delta^2u,u\rangle)\nabla u
+4\nabla(\langle\nabla\Delta u,\nabla u\rangle)\nabla u 
+2(\nabla\langle\nabla\Delta u,u\rangle)\Delta u
+2\Delta(\langle\nabla\Delta u,u\rangle\nabla u)
\\
&
-2\langle\nabla\Delta u,u\rangle\nabla\Delta u
-2\Delta\big(\Delta u|\nabla u|^2\big)+2\nabla\big(|\Delta u|^2\nabla u\big).
\end{align*}
Now, in order to understand when a map of the form \eqref{dfn:tilde-u}
is triharmonic we first of all consider its last component from which
we get the constraint \(\langle \tilde u,P_3(\tilde u)\rangle=0\). For the remaining components of \(\tilde u\) we then get \(P_3(\tilde u)=0\) and we can conclude that maps of the form \eqref{dfn:tilde-u} are triharmonic if \(P_3(\tilde u)=0\) which is precisely the statement of the Lemma.
\end{proof}

With this preparation at hand we can now provide a construction method for triharmonic maps from harmonic maps.

\begin{Satz}
\label{thm:triharmonic-main}
Let $m,\ell\in\mathbb{N}$.
The map
\(q:\mathbb{R}^{m}\setminus\{0\}\rightarrow \s^{m^{\ell}}\)
given by
\begin{align*}
q:=\big(\sin\gamma\cdot u^{(\ell)}_{i_1\dots i_{\ell}},\cos\gamma\big),   
\qquad\gamma\in(0,\pi/2),
\end{align*}
is a proper triharmonic map if and only if the equation
\begin{align}
\label{eq:tri-rotation}
&3\sin^4\gamma\,\ell^3(m+\ell-2)^3-2\sin^2\gamma\,\ell^2(m+\ell-2)^2[4+(m+\ell-6)(4+\ell)+\ell(m+\ell-2)] \\
\nonumber&+\ell(\ell+2)(\ell+4)(m+\ell-2)(m+\ell-4)(m+\ell-6)=0
\end{align}
is satisfied.
\end{Satz}
\begin{proof}
Since the last component of $q$ is constant, it satisfies \eqref{eq:triharmonic-tilde-u} trivially. Below we hence consider the components $\sin\gamma\cdot u^{(\ell)}_{i_1\dots i_{\ell}}$ only, which in the rest of the proof will again be denoted by $q$ for simplicity of notation. 

Straightforward computations give
\begin{align*}
\Delta^3 q=&-\frac{\ell(\ell+2)(\ell+4)}{r^6}(m+\ell-2)(m+\ell-4)(m+\ell-6)q
,\\
\Delta\big(\langle\Delta q,\nabla q\rangle\nabla q\big)=&0
,\\
\nabla\big(\langle\Delta q,\nabla q\rangle\Delta q\big)=&0
,\\
\nabla((\Delta|\nabla q|^2)\nabla q)=&
\sin^2\gamma\frac{2\ell^2(m+\ell-2)^2(m-4)}{r^6}q,
\\
\nabla\big(|\nabla q|^4\nabla q\big)=&-\sin^4\gamma\frac{\ell^3(m+\ell-2)^3}{r^6} q
,\\
(\nabla\langle\Delta^2 q,q\rangle)\nabla q=&
0
,\\
\nabla(\langle\nabla\Delta q,\nabla q\rangle)\nabla q=&
0
,\\
(\nabla\langle\nabla\Delta q,q\rangle)\Delta q=&-\sin^2\gamma\frac{2\ell^2(m+\ell-2)^2(m-4)}{r^6}q
,\\
\Delta(\langle\nabla\Delta q,q\rangle\nabla q)=&0
,\\
\langle\nabla\Delta q,q\rangle\nabla\Delta q=&\sin^2\gamma\frac{4 \ell^2(m+\ell-2)^2}{r^6}q
,\\
\Delta\big(\Delta q|\nabla q|^2\big)=&\sin^2\gamma\frac{(m+\ell-2)^2(m+\ell-6)\ell^2(4+\ell)}{r^6}q
,\\
\nabla\big(|\Delta q|^2\nabla q\big)=&-\sin^2\gamma\frac{\ell^3 (m+\ell-2)^3}{r^6}q.
\end{align*}
Plugging the above results into \eqref{eq:triharmonic-tilde-u} yields the claim. 
\end{proof}

Solving the quadratic equation \eqref{eq:tri-rotation} for \(\sin^2\gamma\) we get
\begin{align}
\label{eq:sinq}
 \sin^2\gamma=\frac{2}{3}+\frac{4(m-5)}{3\ell(m+\ell-2)}\pm\frac{1}{3}\sqrt{1+\frac{2(8-m)}{\ell(\ell+m-2)}-\frac{8(m^2-10m+22)}{\ell^2(\ell+m-2)^2}}.   
\end{align}
For the following considerations it is helpful to study the equations
\begin{align}
\label{eq:sinqm}
 \sin^2\gamma=\frac{2}{3}+\frac{4(m-5)}{3\ell(m+\ell-2)}-\frac{1}{3}\sqrt{1+\frac{2(8-m)}{\ell(\ell+m-2)}-\frac{8(m^2-10m+22)}{\ell^2(\ell+m-2)^2}}   
\end{align}
and
\begin{align}
\label{eq:sinqp}
 \sin^2\gamma=\frac{2}{3}+\frac{4(m-5)}{3\ell(m+\ell-2)}+\frac{1}{3}\sqrt{1+\frac{2(8-m)}{\ell(\ell+m-2)}-\frac{8(m^2-10m+22)}{\ell^2(\ell+m-2)^2}}   
\end{align}
separately.

\smallskip

In the following lemma we discuss under which conditions on $m$ and $\ell$ identity \eqref{eq:sinq} admits solutions.
Note however that, as in the biharmonic case, not every solution of this identity gives rise to a proper triharmonic map: In the construction of Nakauchi \cite{MR4593065}, which we recalled in Section\,\ref{sec:prelim}, we need to impose the condition $\ell\leq m$. In other words, solutions of \eqref{eq:constraint-biharmonic} with $\ell>m$ do not give rise to triharmonic maps.

\begin{Lem}
Equation \eqref{eq:sinq} admits solutions if and only if one of the following conditions is satisfied
\begin{enumerate}
    \item$\ell=1$ and $m\in\{6,7\}$;
    \item $\ell=2$ and $m\in\{3,5,6,7,8,9,10,11\}$; 
    \item $\ell=3$ and  $2\leq m\leq 26$;
    \item $\ell=4$ and $m\geq 3$;
     \item $\ell=5$ and $m\geq 2$;
    \item $\ell\geq 6$ and $m\geq 1$. 
\end{enumerate}
\end{Lem}
\begin{proof}
We study the cases $\ell=1,2,3$ and $\ell\geq 4$ separately.
For $\ell=1,2,3$ we just state the results, they all follow from straightforward computations and considerations.

\smallskip

For $\ell=1$ equation \eqref{eq:sinqm} admits solutions if and only if $m\in\{6,7\}$.
For $\ell=1$ equation \eqref{eq:sinqp} admits no solutions for any $m\in\mathbb{N}$.

\smallskip

For $\ell=2$ equation \eqref{eq:sinqm} has solutions if and only if $m\in\{5,6,7,8,9,10,11\}$. 
For $\ell=2$ equation \eqref{eq:sinqp} has solutions if and only if $m=3$.

\smallskip

For $\ell=3$ equation \eqref{eq:sinqm} admits solutions if and only if $4\leq m\leq 26$.
For $\ell=3$ equation \eqref{eq:sinqp} admits solutions if and only if $m=2,3$.

\smallskip

Let $\ell\geq 4$, $m\geq 8$. Then, we have
\begin{align*}
0<\frac{4(m-5)}{3\ell(m+\ell-2)}\leq\frac{m-5}{3(m+2)}<\frac{1}{3}.    
\end{align*}
Further, we have
\begin{align*}
   0> -\frac{1}{3}\sqrt{1+\frac{2(8-m)}{\ell(\ell+m-2)}-\frac{8(m^2-10m+22)}{\ell^2(\ell+m-2)^2}}\geq -\frac{1}{3}.
\end{align*}
Hence, for $\ell\geq 4$, $m\geq 8$, from \eqref{eq:sinqm} we get
\begin{align*}
 \frac{1}{3}<\sin^2\gamma <1.  
\end{align*}
Therefore, for $\ell\geq 4$, $m\geq 8$ there exists a solution of \eqref{eq:sinqm}.

\smallskip

Straightforward computations using equation \eqref{eq:sinqm} show that for any pair $(\ell,m)$ of natural numbers with $\ell\geq 4$ and $m\in\{1,\dots,7\}$, except for $(\ell,m)\in\{(4,2), (4,1), (5,1)\}$,
there exists $c_1,c_2\in(0,1)$ with $c_1<c_2$ such that $c_1<\sin^2\gamma <c_2$. Thus, for these pairs \eqref{eq:sinqm} admits a solution.

\smallskip

For $(\ell,m)\in\{(4,2), (4,1), (5,1)\}$ equation \eqref{eq:sinqp} admits a solution.
\end{proof}

\begin{Cor}
\label{cor:triharmonic}
For each $m\in\mathbb{N}$ with $m\geq 3$ there exists a natural number $\ell\in\mathbb{N}$ with $\ell\leq m$ such that there is a triharmonic map $\mathbb{R}^{m}\setminus\{0\}\rightarrow \s^{m^{\ell}}$.
\end{Cor}

\section{Polyharmonic maps between spheres from eigenmaps}
\label{sec:eigen}
In this section we provide a construction method for 
proper \(r\)-harmonic maps between spheres by a suitable deformation of eigenmaps.

\medskip

We define the \(r\)-energy for a map between Riemannian manifolds as follows, where we need to distinguish between the even ($r=2s$) and the odd ($r=2s+1$) case:
\begin{align}
\label{eq:r-energy}
E_{2s}(\phi)=&\int_M|\bar\Delta^{s-1}\tau(\phi)|^2\dv, \\
\nonumber E_{2s+1}(\phi)=&\int_M|\bar\nabla\bar\Delta^{s-1}\tau(\phi)|^2\dv.
\end{align}

The Euler-Lagrange equations associated with 
the r-energy functionals \eqref{eq:r-energy} 
have been calculated by Maeta \cite{MR3007953} and
are characterized as follows:

\begin{enumerate}
 \item The critical points of \(E_{2s}\) are given by
\begin{align}
\label{tension-2s}
0=\tau_{2s}(\phi):=&-\bar\Delta^{2s-1}\tau(\phi)-R^N(\bar\Delta^{2s-2}\tau(\phi),d\phi(e_j))d\phi(e_j) \\
\nonumber&+\sum_{l=1}^{s-1}\bigg(R^N(\bar\nabla_{e_j}\bar\Delta^{s+l-2}\tau(\phi),\bar\Delta^{s-l-1}\tau(\phi))d\phi(e_j) \\
\nonumber&\hspace{1cm}-R^N(\bar\Delta^{s+l-2}\tau(\phi),\bar\nabla_{e_j}\bar\Delta^{s-l-1}\tau(\phi))d\phi(e_j)
\bigg).
\end{align}
\item The critical points of \(E_{2s+1}\) are given by
\begin{align}
\label{tension-2s+1}
0=\tau_{2s+1}(\phi):=&\bar\Delta^{2s}\tau(\phi)+R^N(\bar\Delta^{2s-1}\tau(\phi),d\phi(e_j))d\phi(e_j)\\
\nonumber&-\sum_{l=1}^{s-1}\bigg(R^N(\bar\nabla_{e_j}\bar\Delta^{s+l-1}\tau(\phi),\bar\Delta^{s-l-1}\tau(\phi))d\phi(e_j) \\
\nonumber&-R^N(\bar\Delta^{s+l-1}\tau(\phi),\bar\nabla_{e_j}\bar\Delta^{s-l-1}\tau(\phi))d\phi(e_j)
\bigg) \\
&\nonumber-R^N(\bar\nabla_{e_j}\bar\Delta^{s-1}\tau(\phi),\bar\Delta^{s-1}\tau(\phi))d\phi(e_j).
\end{align}
\end{enumerate}
Here, we have set \(\bar\Delta^{-1}=0\), \(\{e_j\},j=1,\ldots,\dim M\) denotes an orthonormal basis of \(TM\)
and we are applying the Einstein summation convention.

\smallskip

It is obvious that due to the large number of terms in \eqref{tension-2s}, \eqref{tension-2s+1} it is a challenging task
to construct explicit solutions of the latter. Hence, it is favorable to choose a different approach in order to construct critical points of \eqref{eq:r-energy} than trying to explicitly solve the Euler-Lagrange equations.

\begin{Lem}
Let \(v\colon\s^m\to\s^{n-1}\) be an eigenmap with energy density
\(|\nabla v|^2=\lambda\).
Then, the map \(q:=\s^m\to\s^n\) defined by
\begin{align}
\label{dfn:q}
q:=(\sin\delta\cdot v,\cos\delta),\qquad \delta\in(0,\frac{\pi}{2})
\end{align}
has \(r\)-energy
\begin{align}
\label{eq:r-energy-q}
E_r(q)=\vol(\s^m,g_R)\lambda^r\epsilon_r,    
\end{align}
where
\begin{align*}
\epsilon_r=\epsilon_r(\delta):=\sin^2\delta\cos^{2(r-1)}\delta.    
\end{align*}
\end{Lem}

Similar recursion formulas have been obtained for polyharmonic hyperspheres,
see \cite[Proposition 2.10]{MR4106647}, \cite[Proposition 3.1]{MR3711937} and in the study of polyharmonic curves in spheres \cite[Lemma 2.2]{MR4542687}.

\begin{proof}
For an arbitrary map \(u\colon M\to\s^n\) 
we can express the higher order derivatives of the tension
field as follows
\begin{align*}
\tau(u)=&\Delta u+|\nabla u|^2u,\\
\bar\nabla\tau(u)=&\nabla\Delta u
+\langle\Delta u,\nabla u\rangle u
+\nabla(|\nabla u|^2u),\\
\bar\Delta\tau(u)=&\Delta^2u+\langle\nabla\Delta u,\nabla u\rangle u
+\nabla(\langle\Delta u,\nabla u\rangle u)
+\Delta(|\nabla u|^2u)+|\nabla u|^4u,\\
\bar\nabla\bar\Delta\tau(\phi)=&\nabla\Delta^2u
+\langle\Delta^2u,\nabla u\rangle u
+\langle \nabla(\langle\nabla\Delta u,\nabla u\rangle u)
+\nabla \tr\nabla_{(\cdot)}\big(\langle\Delta u,\nabla_{(\cdot)}u\rangle u\big) \\
&+\langle \tr\nabla_{(\cdot)}\big(\langle\Delta u,\nabla_{(\cdot)}u\rangle u\big),\nabla u\rangle u
+\nabla\Delta (|\nabla u|^2u) +\langle\Delta(|\nabla u|^2u),\nabla u\rangle u \\
&+\nabla(|\nabla u|^4u).
\end{align*}
If we consider a map of the form \eqref{dfn:q}, then we get
\begin{align*}
\tau(q)=&\lambda\sin\delta\cos\delta\big(-\cos\delta v,\sin\delta\big),\\
\bar\nabla\tau(q)=&
\lambda\sin\delta\cos\delta\big(-\cos\delta\nabla v,0\big),\\
\bar\Delta\tau(q)=&-\lambda\cos^2\delta\tau(q),\\
\bar\nabla\bar\Delta\tau(\phi)=&-\lambda\cos^2\delta\bar\nabla\tau(\phi).
\end{align*}
By iteration we thus find 
\begin{align*}
\bar\Delta^k\tau(q)&=(-1)^k\lambda^k\cos^{2k}\delta\tau(q),\\
\bar\nabla\Delta^k\tau(q)&=(-1)^k\lambda^k\cos^{2k}\delta \bar\nabla\tau(q).
\end{align*}
Inserting these identities into \eqref{eq:r-energy-q} already completes the proof.
\end{proof}

In the following theorem we provide a construction method for 
proper \(r\)-harmonic maps between spheres from eigenmaps.

\begin{Satz}
\label{thm:r-harmonic-main}
Let \(v\colon\s^m\to\s^{n-1}\) be an eigenmap.
Then, the map \(q:\s^m\to\s^n\) defined by
\begin{align*}
q:=(\sin\delta\cdot v,\cos\delta)
\end{align*}
is proper \(r\)-harmonic if and only if \(\sin\delta=\frac{1}{\sqrt{r}}\). These maps are unstable critical points of the \(r\)-energies \eqref{eq:r-energy-q}.
\end{Satz}
\begin{proof}
Using the reduced expression for the \(r\)-energy provided by \eqref{eq:r-energy-q} it is now straightforward to calculate the critical points of the \(r\)-energy.
Now, a direct calculation shows that
\begin{align*}
\frac{d}{d\delta}\epsilon_r(\delta)=2\sin\delta\cos^{2r-3}\delta(1-r\sin^2\delta)    
\end{align*}
and this expression clearly vanishes if \(\sin\delta=\frac{1}{\sqrt{r}}\) leading to the first claim.

In order to show that these critical points are actually unstable
we compute
\begin{align*}
\frac{d^2}{d\delta^2}\epsilon_r(\delta)  
=\big(2\cos^{2r-2}\delta-2(2r-3)\sin^2\delta\cos^{2r-4}\delta\big)
(1-r\sin^2\delta)
-4\sin^2\delta\cos^{2r-2}\delta
\end{align*}
and it is easy to see that
\begin{align*}
\frac{d^2}{d\delta^2}\epsilon_r(\delta)\big|_{\sin\delta=\frac{1}{\sqrt{r}}}<0     
\end{align*}
leading to the instability statement and thus completing the proof.
\end{proof}

\begin{Bem}
Theorem\,\ref{thm:r-harmonic-main} has been known for $r=2$, see for example \cite{bra25}.  
\end{Bem}

\bibliographystyle{plain}
\bibliography{mybib}
\end{document}